\documentclass{article}
\usepackage{spconf,amsmath,graphicx,amssymb}
\usepackage{float}

\title{Asynchronous Algorithms for Solving Linear Programs}

\name{Tarek A. Lahlou \quad Thomas A. Baran \thanks{The authors wish to thank Analog Devices, Bose Corporation, and Texas Instruments for their support of innovative research at MIT and within the Digital Signal Processing Group.}}
\address{Digital Signal Processing Group \\Massachusetts Institute of Technology }

\def\a{\underline{a}}
\def\b{\underline{b}}
\def\c{\underline{c}}
\def\d{\underline{d}}
\def\e{\underline{e}}
\def\f{\underline{f}}
\def\g{\underline{g}}
\def\h{\underline{h}}
\def\k{\underline{k}}
\def\xopt{\underline{x}^{\star}}
\def\w{\underline{w}}
\def\x{\underline{x}}
\def\y{\underline{y}}
\def\z{\underline{z}}

\DeclareMathOperator*{\minimize}{minimize}

\begin{document}
\ninept
\maketitle

\begin{abstract}
	In this paper we design and analyze algorithms for asynchronously solving linear programs using nonlinear signal processing structures. In particular, we discuss a general procedure for generating these structures such that a fixed-point of the structure is within a change of basis the minimizer of an associated linear program. We discuss methods for organizing the computation into distributed implementations and provide a treatment of convergence.  The presented algorithms are accompanied by numerical simulations of the Chebyshev center and basis pursuit problems.
\end{abstract}


\section{Introduction} 
\label{sec:introduction}
	Traditional linear programming algorithms typically iterate an ordered sequence of operations until either a suitable stopping criterion is met or a minimizer is identified, e.g.~interior point and basis exchange methods \cite{BertsimasTsitsiklisLP}. In transforming such an algorithm into the distributed setting, issues related to task allocation, communication, and global versus local synchronization become of paramount importance. Computational reorganization often involves partitioning the sequence of operations into those which distribute in the sense that their execution can make use of many processors that do not interchange information and those which do not, 
	e.g. methods like \cite{NealBoyd}-\cite{WeiOzdaglar1} and the distributed simplex method in \cite{DistributedSimplex}.
	Global optimization algorithms in the style of monte carlo methods paired with local search techniques make efficient use of distributed resources but often require careful parameter tuning for a given problem \cite{GlobalAsynchOpt}. Understanding the complexity associated with a sequential algorithm differs in many ways from its distributed counterpart as the penalties incurred for communication and synchronization may be more expensive than that suggested by traditional complexity measures \cite{ParallelComplexity}. 

	In this paper we specifically address linear programs from a conservative signal processing perspective consistent with the general framework presented in \cite{BaranLahlouPartI}. In particular, we design and analyze a generally nonlinear signal-flow structure which when implemented asynchronously solves an associated system of stationarity conditions that was developed in \cite{BaranThesis}. The viewpoint of solving either an optimization or constraint satisfaction problem using an asynchronous signal processing system naturally lends itself to understanding important issues such as algorithmic scalability, robustness with respect to communication and processing delays, and computational heterogeneity. Furthermore, issues pertaining to the identification of sufficient conditions for convergence may be addressed using well-known stability and dynamic systems results.

\section{Preliminaries}
\label{sec:preliminaries}
	The optimization problem of minimizing a linear objective function subject to linear equality and inequality constraints is described in {\it standard form} as
	\begin{equation}  \label{eq:standard-form}
		\begin{array}{rl} \displaystyle
			 \minimize_{\x} & \f^T\x \\ 
			\text{subject to} & A\x \leq \b \\ 
			& \,\,\,\,\x \geq \underline{0}
		\end{array}
	\end{equation}
	where $\x\in\mathbb{R}^N$ is the decision vector, $\f\in\mathbb{R}^N$ is the cost vector, $\b\in\mathbb{R}^M$ is the constraint vector, and $A\in\mathbb{R}^{M\times N}$ is the coefficient matrix. We describe an asynchronous algorithm for solving \eqref{eq:standard-form} in this paper via a signal-flow structure consisting of subsystems realized as maps coupled together using asynchronous delays, i.e.~randomly triggered sample-and-hold elements which we model as independent Bernoulli processes. The general strategy underlying the presented algorithms described in this way is to determine a solution to a system of stationarity conditions associated with a particular reformulation of \eqref{eq:standard-form} by interconnecting linear and nonlinear signal-flow elements such that they collectively describe the behavior of the stationarity conditions, resolving any delay-free loops, and finally running the system to a fixed-point. The delay-free loops are either resolved algebraically, e.g. using the automated techniques in \cite{BaranLahlouPartICASSP}, or by inserting asynchronous delays depending on the specific form of the associated nonlinearity.	 

	We refer to a linear program problem statement organized according to the following conventions as being in {\it asynchronous form}: 
	\begin{enumerate}
		\item[(i)] 		every vector (except the cost vector) is involved in a system of linear equations,
		\item[(ii)] 	every vector in (i) is either fixed, unconstrained, or non-negative,
		\item[(iii)] 	every vector in (i) with non-zero cost coefficients must be unconstrained.
	\end{enumerate}
	Indeed, a linear program in standard form can always be recast into asynchronous form as
	\begin{equation}  \label{eq:asynchronous-form}
		\begin{array}{rl} \displaystyle
			\minimize & \w^T\z_1 \\ 
			\text{subject to} & B\z_1 = \z_2 \\ 
			& \,\,\,\,\z_2 \geq \underline{0}
		\end{array}
	\end{equation}
	where the minimization is explicitly over those variables which are unconstrained and non-negative and 
	\begin{equation*}
		B 			=  	\left[	\begin{array}{cc} 	0 & I \\ I & -A 		\end{array}\right], \hspace{.025in} 
		\w			= 	\left[	\begin{array}{c}	\underline{0} \\ \f 	\end{array}\right], \hspace{.025in} 
		\z_{1} 		= 	\left[	\begin{array}{c}	\b   \\ \x_1			\end{array}\right], \hspace{.025in} 
		\z_{2} 		= 	\left[	\begin{array}{c}	\x_2 \\ \y 				\end{array}\right].
	\end{equation*}
	Note that we have introduced a non-negative vector $\y\in\mathbb{R}^{M}$ in order to enforce the linear inequality $A\x \leq \b$ and have made two equality-constrained copies of $\x$ where $\x_1$ encodes the non-zero cost coefficients and $\x_2$ enforces the non-negativity constraints.
	
	Consistent with the presentation in \cite{BaranLahlouPartI}\cite{BaranLahlouPartII}, the system of stationarity conditions for the formulation in \eqref{eq:asynchronous-form} is of the general form 
	\begin{eqnarray}
		\d^{\star}  & = &  G\c^{\star}         \label{eq:stationarity-linear}		\\ 
		\c^{\star}  & = &  m(\d^{\star})	   \label{eq:stationarity-nonlinear} 
	\end{eqnarray}
	where $G$ is an orthogonal matrix, $m(\cdot)$ is an element-wise memoryless nonlinearity, $\c$ and $\d$ denote respectively the input to and output from $G$, and $\c^{\star}$ and $\d^{\star}$ denote a fixed-point of the algebraic system.  The vectors $\c$ and $\d$ are in particular organized according to the ordering of the inputs followed by the outputs of the linear equality constraints in \eqref{eq:asynchronous-form}, i.e. 
	\begin{equation*}
	\d = \left[\begin{array}{c}\d_{\z_{1}} \\ \d_{\z_{2}}\end{array}\right] \quad \text{and} \quad \c = \left[\begin{array}{c}\c_{\z_{1}}\\\c_{\z_{2}}\end{array}\right].
	\end{equation*}
	Furthermore, the matrix $G$ in \eqref{eq:stationarity-linear} is generated using $B$ in \eqref{eq:asynchronous-form} as 
	\begin{eqnarray}
		G &=& \left(I + R\right) \left( I - R\right)^{-1}     \label{eq:G}  \\
		  & = & 	\left(I+R\right)^{2}\left[\begin{array}{cc}
					\left(I+B^{T}B\right)^{-1} & 0\\
					0 & \left(I+BB^{T}\right)^{-1}
					\end{array}\right] \nonumber
	\end{eqnarray}
	where $R=-R^T$ is a skew-symmetric matrix of the form
	\begin{equation*}
		R = \left[\begin{array}{cc} 0 & -B^{T}\\ B & 0 \end{array}\right]. 
	\end{equation*}
	The orthogonality of $G$ is readily verified using the fact that $I+R$ and $I-R$ commute. Moreover, it follows that $G$ is an element of the subset of special orthogonal matrices which do not have eigenvalues of $-1$. When $B$ is itself an orthogonal matrix, e.g. the discrete Fourier transform matrix, the expression in \eqref{eq:G} simplifies to $G=R$ and thus if a fast implementation of $B$ is available it may be used in the implementation of $G$.

	The memoryless nonlinearities used to realize the stationarity conditions in \eqref{eq:stationarity-nonlinear} for a linear program described in asynchronous form are listed in Table~1 for a single element $z$ which either belongs to $\z_1$ (input) or $\z_2$ (output). The nonlinearities are specifically described for the various pairings of cost contributions and set memberships which arise in linear programs cast this way.  For example, $z$ being unconstrained with no contribution to the overall cost follows from the second row with $\rho=0$. In the sequel we resolve delay-free loops in the presented signal-flow structures associated with the first two rows algebraically and the third row using asynchronous delays.
	\begin{center}
		\textsc{Table 1. nonlinearities for linear programs}\vspace{.05in}
		\begin{tabular}{ l | l | l | l }
			cost     & set membership & $m(\cdot)$ (input) & $m(\cdot)$ (output)\\
	  		\hline \hline
	  		none     & 	fixed $z = \rho$	& 	$c = -d + 2\rho$ 	& 	$c = d - 2\rho$  	\\
	  		$\rho z$ & 	unconstrained       & 	$c = d-2\rho$    	&	$c=-d+2\rho$ 	   	\\
	  		none     & 	non-negative 		& 	$c=|d|$      		& 	$c=-|d|$ 	  	   	\\
		\end{tabular}
	\end{center} 
	
	Given a fixed-point $\c^{\star}$ and $\d^{\star}$ of the stationarity conditions \eqref{eq:stationarity-linear}-\eqref{eq:stationarity-nonlinear}, the argument $\xopt$ which minimizes \eqref{eq:asynchronous-form} and consequntly \eqref{eq:standard-form} may be determined by extracting either $\x_1$ or $\x_2$ from 
	\begin{eqnarray}\label{eq:stationary-solution}
		\z_{1}^{\star}  = \frac{1}{2}\left(\d_{\z_1}^{\star} + \c_{\z_1}^{\star}\right) \quad \text{ or } \quad
		\z_{2}^{\star}  = \frac{1}{2}\left(\d_{\z_2}^{\star} - \c_{\z_2}^{\star}\right). 
	\end{eqnarray}

\section{Asynchronous Linear Program Algorithms}
\label{sec:a-general-asynchronous-algorithm}
\begin{figure}[t]
  		\centering
  		\centerline{\includegraphics[width=8cm]{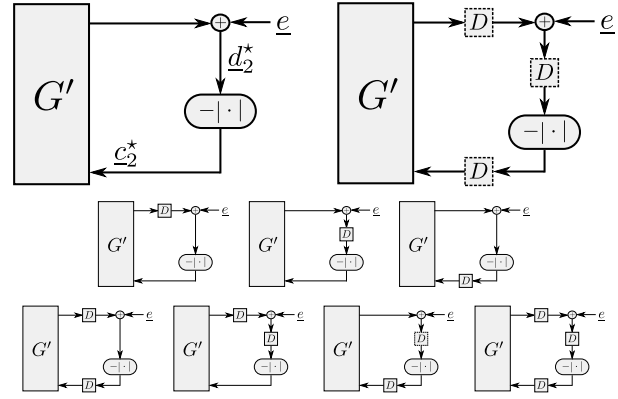}}\caption{A general signal-flow structure (top left) used in the initial description of the stationarity conditions for linear programs in asynchronous form and the resulting signal-flow structure (top right) indicating three possible locations for distributing system state specifically indicated using dashed boxes.}
		\label{fig:general-signal-flow-structures}\vspace{-0.08in}
	\end{figure}

	In this section we present a general asynchronous algorithm in the form of a signal-flow structure for solving linear programs described in asynchronous form consistent with the general strategy previously described. We begin this presentation by partitioning the linear stationarity conditions in \eqref{eq:stationarity-linear} where we specifically delineate between those variables associated with memoryless nonlinearities in \eqref{eq:stationarity-nonlinear} which are affine maps and those which are generally not, i.e.
	\begin{equation}
		\left[\begin{array}{c} \d^{\star}_1 \\ \d^{\star}_2 \end{array}\right] = \left[\begin{array}{cc}G_{11} & G_{12} \\ G_{21} & G_{22} \end{array}\right]\left[\begin{array}{c}\c^{\star}_1 \\ \c^{\star}_2 \end{array}\right]\label{eq:dynamic-system}
	\end{equation}
	where the variables are ordered according to 
	\begin{eqnarray*}\label{eq:c-d-partitioning}
		\d_1  =  \left[	\begin{array}{c} 	\d_{b}    	\\  	\d_{\x_1}   \end{array}		\right]\hspace{-0.025in}, 
		\d_2  =  \left[	\begin{array}{c} 	\d_{\x_2} 	\\  	\d_{\y}     \end{array}		\right]\hspace{-0.025in}, 
		\c_1  =  \left[	\begin{array}{c} 	\c_{\b}   	\\  	\c_{\x_1}	\end{array}		\right]\hspace{-0.025in}, 
		\c_2  =  \left[	\begin{array}{c} 	\c_{\x_2} 	\\  	\c_{\y}		\end{array}		\right]\hspace{-0.025in},
	\end{eqnarray*}
	and where $G$ is generated using \eqref{eq:G} and block partitioned accordingly. Eliminating those variables associated with $\b$ and $\x_{1}$ in \eqref{eq:stationarity-linear}, i.e. $\d_1$ and $\c_1$, results in an affine system of the form
	\begin{eqnarray}
		\d^{\star}_2 & = & G'\c^{\star}_2  + \e \label{eq:asynch-linear}
	\end{eqnarray}
	where
	\begin{eqnarray}
		G' & = &  G_{22}+G_{21}\left(I-SG_{11}\right)^{-1}SG_{12}  \label{eq:g-prime}\\
		\e & = & 2G_{21}\left(I-SG_{11}\right)^{-1}\left[\begin{array}{c}\b \\ -\f\end{array}\right] \\
		S  & = & \left[\begin{array}{cc} -I_{M} & 0 \\ 0 & I_{N}\end{array}\right].
	\end{eqnarray}
	The stationarity conditions corresponding to the described reduced representation are given by \eqref{eq:asynch-linear} paired with  $\c_2^{\star}  =  -|\d_2^{\star}|$. A signal-flow structure depicting a fixed-point of these stationarity conditions is portrayed in Figure~\ref{fig:general-signal-flow-structures} on the top left. Another class of related algorithms follow from defining a sequence of signal-flow systems which smoothly deform into this structure. Referring again to Fig.~\ref{fig:general-signal-flow-structures}, the signal-flow structure on the top right depicts three possible locations at which asynchronous delays may be inserted in order to break delay-free loops where the dashed boxes labeled $D$ denote a vector asynchronous delay element. This structure forms the basis from which various implementations may be synthesized. Seven natural organizations of the system state are additionally depicted.

	\begin{figure}[t]
		\begin{centering}
		\includegraphics[width=3.2in]{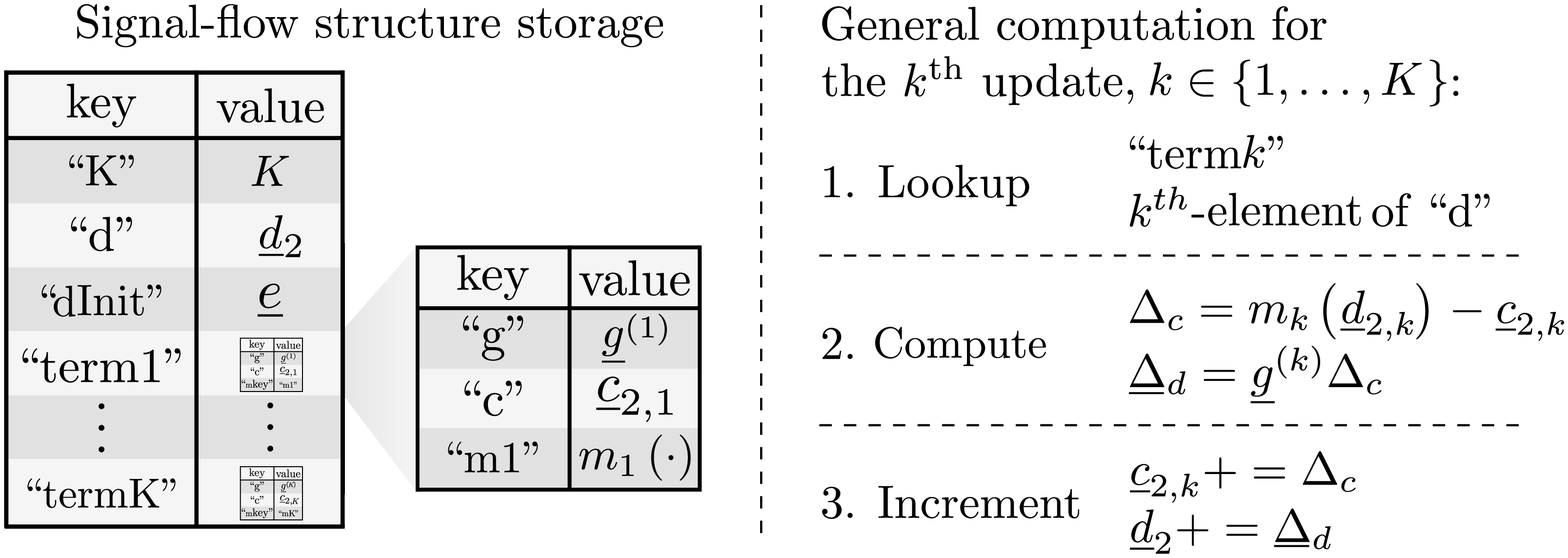}
		\par
		\end{centering}\vspace{-0.15in}
		\caption{Left: a signal processing system organized into a conceptual associative array. Right: the general computation for the $k^{\text{th}}$ update.\label{fig:associative-array}}
		\vspace{-.15in}
	\end{figure}

\section{Distributed asynchronous implementations}\vspace{-0.08in}
\label{sec:asynchronous-implementations}
	We next briefly present an example implementation of an asynchronous signal-flow structure using a conceptual associative array organized using the key-value pairs depicted on the left of Fig.~\ref{fig:associative-array}. The algorithms used in Section~\ref{sec:numerical-simulations} were specifically implemented using this approach where the general computation associated with an asynchronous update is depicted on the right and explained below.

	Consider the causal system, i.e.~with initial condition $\c_{2}[0]=\underline{0}$, defined by the recurrence relation
	\begin{eqnarray} \label{eq:implementation1}
		\begin{array}{rlc}
			\c_2[n] 	&	=	& 		m(\d_2[n-1]) \\
			\d_{2}[n] 	&	=	& 		G'\c_{2}[n] + \e 
		\end{array}
	\end{eqnarray}
	and note that a fixed-point of this system is a fixed-point of \eqref{eq:asynch-linear}. For $n\geq1$ an equivalent description obtained via manipulation is\vspace{-0.04in}
	\begin{equation}
	\d_{2}[n] = \d_2[n-1] +\sum_{k=1}^{K} \g^{(k)}\left(m(\d_{2,k}[n-1])-\c_{2,k}[n-1]\right) \hspace{-0.002in} \label{eq:implementation2}
	\end{equation}
	\vspace{-0.08in}

	\noindent where $\g^{(k)}$ denotes the $k^{\text{th}}$ column of $G'$ and $(\cdot)_{2,k}$ denotes the $k^{\text{th}}$ element of $(\cdot)_2$. The initial condition $\d_{2}[0]=\e$ is required for \eqref{eq:implementation1} and \eqref{eq:implementation2} to produce the same output $\d_{2}[n]$ for all $n\geq0$. An asynchronous implementation of this system then follows from computational nodes executing the general computation described in Fig.~\ref{fig:associative-array} on the right for a randomly selected $k$. The ``compute'' stage computes the $k^{\text{th}}$ term of the summation in \eqref{eq:implementation2} using information obtained in the ``lookup'' stage, while the ``increment'' stage is used to update $\d_2$ and $\c_{2,k}$ without requiring a full read-and-write operation.

\vspace{-0.07in}
\section{Convergence analysis}
\begin{figure}[t]
  		\centering
  		\centerline{\includegraphics[width=2.4in]{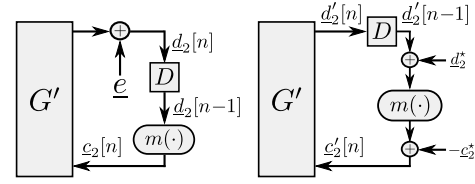}}\vspace{-0.12in}\caption{The signal-flow structure used for convergence analysis.}
		\label{fig:convergence}\vspace{-0.18in}
	\end{figure}

\vspace{-0.07in}
\label{sec:convergence analysis}
	In this section we analyze convergence properties of the signal-flow structure in Fig.~\ref{fig:convergence} on the left for both synchronous and asynchronous settings. To facilitate this analysis, let $T$ denote an iterated operator mapping the output of the delay to the input such that composing $T$ with itself  {\it ad infinitum} produces a synchronous implementation, i.e.
	\begin{eqnarray}
		T(\d_{2}) = G'm(\d_2) + \e. \label{eq:T}
	\end{eqnarray}
	Indeed, for the given operator linear convergence to a unique fixed-point is guaranteed provided that $T$ is Lipschitz continuous with constant $L_T\in[0,1)$, i.e. for every $\d_2^{(1)}$ and $\d_{2}^{(2)}$ $T$ satisfies\vspace{-0.01in}
	\begin{eqnarray}
		\left\|T\left(\d_{2}^{(1)}\right)- T\left(\d_{2}^{(2)}\right) \right\|_2 \leq L_T \left\|\d_{2}^{(1)} - \d_{2}^{(2)} \right\|_2.\label{eq:contractive-inequality}
	\end{eqnarray}
	\vspace{-0.11in}

	\noindent
	Establishing the existance of a fixed-point for this case follows immediately by assigning $\d_{2}^{(1)}=\d_{2}[n]$ and  $\d_{2}^{(2)} = \d_{2}[n-1]$ and taking the limit of the iterated inequality, i.e.~
	\begin{eqnarray*}
			\lim_{n\rightarrow\infty} \left\Vert \d_{2}\left[n\right]-\d_{2}\left[n-1\right]\right\Vert_2 \leq \left\Vert \d_{2}\left[1\right]-\d_{2}\left[0\right]\right\Vert_2\lim_{n\rightarrow\infty}L_T^{n-1} = 0. \label{eq:Lipschitz}
	\end{eqnarray*}
	Suppose two signals $\underline{d}_{2}^{(1)}[n]\neq \underline{d}_{2}^{(2)}[n]$ are both fixed-points of $T$. Taking the limit of the iterated inequality in \eqref{eq:contractive-inequality} for this case yields
	\begin{equation*}
		\lim_{n\rightarrow\infty}\left\Vert \d^{(1)}_{2}\left[n\right]-\d^{(2)}_{2}\left[n\right]\right\Vert_2 \leq\left\Vert \d^{(1)}_{2}\left[0\right]-\d^{(2)}_{2}\left[0\right]\right\Vert_2\lim_{n\rightarrow\infty}L_T^{n} = 0.\label{eq:unique-ineq}
	\end{equation*}
	Since a fixed-point is defined as being independent of $n$ we have $\d^{(1)}_{2}[n]=\d^{(2)}_{2}[n]$ which contradicts the original assumption, thereby establishing that the fixed-point of $T$ is unique.

	We proceed to analyze the asynchronous setting using the system in Fig.~\ref{fig:convergence} on the right consisting of the difference between the system on the left and a fixed-point. Note that when $m(\cdot)$ satisfies \eqref{eq:Lipschitz} with Lipschitz constant $L_m$ then the map $m'$ from $\d_2'[n-1]$ to $\c_2'[n]$ does too. Let $A_{n-m}$ denote the event that a particular asynchronous delay element last fired at time $n-m$ and let $h[m]$ denote the corresponding geometric distribution with parameter $p\in(0,1)$. Then, using the law of total expectation on $\c'_{2,k}[n]$ yields
	\begin{eqnarray*} 
	\begin{array}{lll} 
		E\left[\left(\c_{2,k}'[n]\right)^{2}\right] 	&	=	& 	\displaystyle\sum_{m=0}^{\infty}	h[m]	E\left[\left(\c_{2,k}'[n]\right)^{2} \mid A_{n-m}\right] \\ 
								& 	=	& 	\displaystyle\sum_{m=0}^{\infty}	h[m]	E\left[\left(m'_{k}(\d'_{2,k}[n-m-1])\right)^{2}\right].
	\end{array}
	\end{eqnarray*}
	Summing the $K$ elements composing $\c_{2}$, interchanging the summations, and rewriting the expression using vectors yields
	\begin{eqnarray}
		\displaystyle E\left[\left\|\c'_{2}[n]\right|^2_2\right] 	&=& \displaystyle \sum_{m=0}^{\infty}h[m]	E\left[\left\|m'\left(\d'_{2}[n-m-1]\right)\right\|_2^2\right]\label{eq:lipschitz-pre-inequality} \\ 
								& 	\leq	& 	\displaystyle L_{m}^{2} 	\sum_{m=0}^{\infty}	h[m]	E\left[\left\|\d'_{2}[n-m-1]\right\|_2^2\right] \label{eq:lipschitz-post-inequality}
	\end{eqnarray}
	where we have explicitly used the Lipschitz continuity of $m'$ in going from \eqref{eq:lipschitz-pre-inequality} to \eqref{eq:lipschitz-post-inequality}. Finally, rearranging indices and noting that $E\left[\|\d_2'[n]\|_2^2\right]=E\left[\|\c'_{2}[n]\|_2^2\right]$ leads to the expression 
	\begin{eqnarray*}
		E\left[\left\|\d'_{2}[n]\right\|^2_2\right]	& 	\leq	& 	L_{m}^{2}	\sum_{m=1}^{\infty}	h[m-1]	E\left[\left\|\d'_{2}[n-m]\right\|_2^2\right]
	\end{eqnarray*}
	which takes the form of a convolution bound and is sufficient for asymptotic average convergence if $L_{m} \in[0,1)$.

	For linear programs, however, $T$ is readily verified to satisfy $L_T=1$ since $G'$ is orthogonal and $m(\cdot)$ is non-expansive. Justifying the orthogonality of $G'$ in \eqref{eq:g-prime} follows from showing that $G'$ is a linear isometry. In particular, since $G$ and $S$ are orthogonal matrices
	\begin{eqnarray*}
		 \left\Vert \d^{(1)}_{1}-\d^{(2)}_{1}\right\Vert_2 ^{2} + \left\Vert \underline{d}^{(1)}_{2}-\underline{d}^{(2)}_{2}\right\Vert_2 ^{2}	= 	\left\Vert \underline{c}^{(1)}_{1}-\underline{c}^{(2)}_{1}\right\Vert_2 ^{2}+\left\Vert \underline{c}^{(1)}_{2}-\underline{c}^{(2)}_{2}\right\Vert_2 ^{2}
	\end{eqnarray*}
	and since $\c_1 = S\d_1 + \k$ for some constant $\k$ it follows that
	\begin{equation*}
		\left\Vert \d^{(1)}_{2}-\d^{(2)}_{2}\right\Vert_2 = \left\Vert \c^{(1)}_{2}-\c^{(2)}_{2}\right\Vert_2.
	\end{equation*} 
	Since the operator $T$ as described is on the boundary of guaranteed convergence, it stands to reason that a homotopic relaxation may remedy divergent or oscillatory behavior. Toward this end, define an operator $T'$ with homotopy parameter $\alpha\in[0,1]$ as
	\begin{equation*}
		T'(\d_2,\alpha) = \alpha T(\d_2) + (1-\alpha)T_0(\d_2)
	\end{equation*} 
	where $T_0$ is chosen such that $L_{T_0}\in[0,1)$ and varying $\alpha$ from $0$ to $1$ corresponds to a smooth deformation from $T_0$ to $T$. Therefore, it follows that synchronous or asynchronous implementations of $T'$ as $\alpha\rightarrow 1$ converge to a fixed-point.

\cleardoublepage

\section{Numerical simulations}
\label{sec:numerical-simulations}
	In this section the convergence properties associated with asynchronously solving two linear programming problems using the algorithms developed in this paper are explored where the asynchronous delays were numerically simulated using discrete-time sample-and-hold elements triggered by independent Bernoulli processes. For the sake of comparison, we use the metric of an equivalent iteration to normalize between various probabilities of sampling, i.e.~an equivalent iteration is the same total amount of computation associated with a synchronous iteration where all delays fire independent of the probability of sampling. 
	\begin{figure}[t]
  		\centering
  		\centerline{\includegraphics[width=9cm]{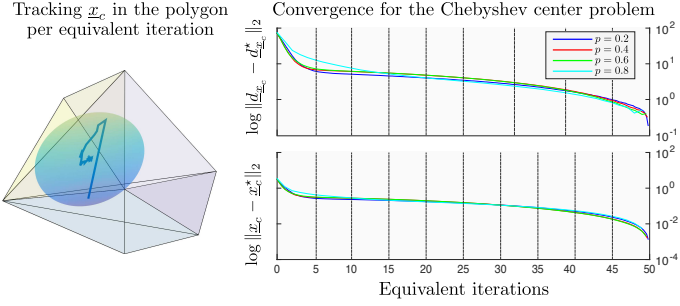}}
  		\vspace{-0.15in}
  		\caption{Asynchronous convergencefor the Chebyshev center problem averaged over $500$ trials.}
  		\vspace{-0.125in}
		\label{fig:cheby-convergence}
	\end{figure}

	\subsection{The Chebyshev center problem}
		Consider as an example the Chebyshev center problem \cite{BoydCVXBook} given by
		\begin{eqnarray} \label{eq:chebyshev}
			\begin{array}{rll} \displaystyle
				\minimize_{\x_c, \, r} 	& 	-r  &\\
				\text{subject to} & \a_i^T\x_c + \|\a_i\|_2r \leq b_i & 1\leq i \leq M \\
				& r \geq 0 & 
			\end{array}
		\end{eqnarray}
		In this form \eqref{eq:chebyshev} is explicitly identifying the largest Euclidean ball which can be inscribed within a convex polytope described in half-space representation by $\{\z \colon A\z\leq b\}$ where $r$ and $\x_c$ denote the balls radius and center, respectively, and $\a_i$ is the $i^{\text{th}}$ column of $A$.  We proceed by recasting this problem into asynchronous form as
		\begin{eqnarray}\label{eq:chebyshev-center}
			\begin{array}{rl} \displaystyle
				\minimize_{\x_c, \z, r_1, r_2} & -r_1 \\
				\text{subject to} &  \left[\begin{array}{ccc}1 & 0 & 0 \\ \underline{n} & -A & b \end{array}\right]\left[\begin{array}{c}r_1 \\ \x_c \\ t \end{array}\right] = \left[\begin{array}{c}r_2 \\ \z\end{array}\right] \\
				& r_2 \geq 0, \,\,\,\, \z \geq \underline{0},\,\,\,\, t=1
				
			\end{array}
		\end{eqnarray}
		where the $i^{\text{th}}$ entry of $\underline{n}$ is $\|\a_i\|_2$. The organization of variables for this problem along with the memoryless nonlinearities associated with the stationarity conditions in \eqref{eq:stationarity-nonlinear} is summarized in Table~2. The matrix in $\eqref{eq:chebyshev-center}$ is used to generate $G$ using \eqref{eq:G} and in the process of synthesizing an asynchronous algorithm adhering to the presented framework the system variables associated with $r_1$, $\x_c$, and $t$ are eliminated using \eqref{eq:asynch-linear} for an appropriate choice of $\d_1$, $\c_1$, $\d_2$, and $\c_2$. Continuing with this notation, the balls center $\x_c^{\star}$ may be recovered from $\c_2^{\star}$ by first recovering $\c_1^{\star}$ and $d_1^{\star}$ via
		\begin{eqnarray}
			\d_1^{\star} & = & \left(I-G_{11}S\right)^{-1}\left(G_{12}\c_2^{\star} + G_{11}\h\right) \\
			\c_1^{\star} & = & S\d_1^{\star} + \h
		\end{eqnarray} 
		followed by partitioning $\z_1^{\star}$ in \eqref{eq:stationary-solution} where
		\begin{eqnarray}
			\h = \left[\begin{array}{c}2\\\underline{0}\\2\end{array}\right]\hspace{1em}\text{and}\hspace{1em}S = \left[\begin{array}{cc}I & 0 \\ 0 & -1\end{array}\right].
		\end{eqnarray}
		
		Convergence results for the Chebyshev center problem in \eqref{eq:chebyshev-center} averaged over $500$ runs of an asynchronous algorithm are depicted in Fig.~\ref{fig:cheby-convergence} on the right where the convex polygon is defined using $200$ random hyperplanes in a $100$-dimensional space. As a further example, the geometric figure on the left tracks the center of the Euclidean sphere in $3$-dimensions for the depicted polygon over the course of an asynchronous algorithm as it converges to its fixed-point. The center is tracked using a blue line and the final Euclidean sphere is also depicted.  
		\begin{center}
			\textsc{Table 2. variable organization for \eqref{eq:chebyshev-center}}\vspace{.05in}
			\begin{tabular}{ l | l | l | l | l }
				variable & cost & set membership & I/O  & implementation \\
			  	\hline \hline                      
			  	$r_{1}$ 	& 	$-r_1$ 	& 	unconstrained 	& 	input  	& 	$ c_{r_1} = d_{r_1} + 2$ 	\\
			  	$\x_c$ 		& 	none 	& 	fixed        	& 	input   	& 	$c_{x} = \d_{x}$			\\
			  	$t$ 		& 	none  	& 	fixed 		  	& 	input 	& 	$c_{t} = -d_{t} + 2$ 		\\
			  	$r_{2}$ 	& 	none  & 	non-negative  	& 	output 	& 	$ c_{r_2} = -|d_{r_2}|$ 	\\
			  	$\z$    	& 	none  & 	non-negative 	& 	output  	& 	$\c_{z} \,\,= -|\d_{z}|$	\\
			\end{tabular}
		\end{center}

	\begin{figure}[t]
	  		\centering
	  		\centerline{\includegraphics[width=7.5cm]{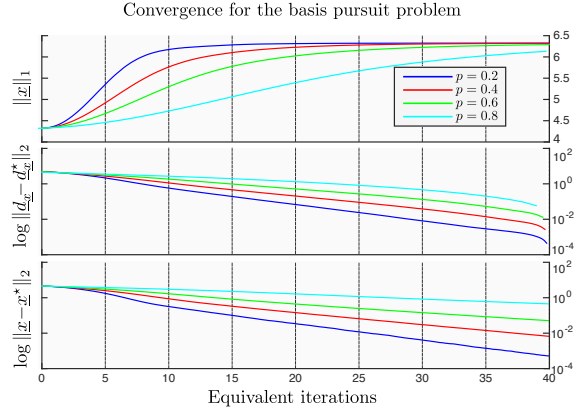}}
	  		\vspace{-0.125in}
	  		\caption{Asynchronous convergence for the basis pursuit problem averaged over $1000$ trials.}
	  		\vspace{-0.17in}
			\label{fig:bp-convergence}
		\end{figure}
	\subsection{The basis pursuit problem}
		Consider as another example the basis pursuit problem given by
		\begin{eqnarray} 
			\begin{array}{rl} \displaystyle
				\minimize_{\x} 		& 	\|\x\|_1 	\\
				\text{subject to}	& 	A\x = \b
			\end{array} 
		\end{eqnarray}	
		Although the process of recasting this problem as a linear program in standard form is well known, we proceed by using a nonlinearity which corresponds to an unconstrained variable $z$ with cost contribution $|z|$. In particular, for $z$ being an input to the problems linear equality constraints the associated nonlinearity is
		\begin{equation*}
			c = m_{1}(d) = 
			\left\{ \begin{array}{rl}
				d+2, 	& 	d<-1 		\\
				 -d, 	& 	|d|\leq1 	\\
				d-2, 	& 	d>1
			\end{array}\right.
		\end{equation*}
		and for $z$ being an output the nonlinearity is $c=-m_{1}(d)$. Convergence results for a $16$-sparse $\xopt$ in a $512$-dimensional space recovered using $200$ random measurements averaged over $1000$ runs of an asynchronous algorithm adhering to the presented framework are depicted in Fig.~\ref{fig:bp-convergence}. In particular, the objective value, log of the distance between $\d_{\x}$ to $\d_{\x}^{\star}$, and log of the distance between $\x$ to $\xopt$ are depicted as a function of equivalent iteration where the delays fire with probability $p=0.2, 0.4, 0.6$, and $0.8$. A homotopic method was used to encourage initial convergence, i.e. the nonlinearity used for equivalent iterations $k=1,\dots,10$ was $c=(1-0.95^{k^{2}})m_1(d)$ and $c=m_1(d)$ thereafter.

\bibliographystyle{IEEEbib}
\bibliography{refs}

\end{document}